\documentclass[fleqn,reqno]{amsart}%{gCOM2e}
\usepackage{mathtools,fancyvrb,breqn}
\usepackage{url}
\usepackage{fullpage}
\fvset{xleftmargin=2\parindent}

\theoremstyle{plain}

\theoremstyle{remark}

\newcommand*{\arraystrut}{\rule[-2mm]{0mm}{6mm}}
\newcommand*{\N}{\noindent}

\begin{document}
\allowdisplaybreaks

%\jvol{00} \jnum{00} \jyear{2013} \jmonth{February}

%\articletype{GUIDE}

\title{Developing explicit Runge-Kutta formulas using open-source software}

\author{Alasdair McAndrew}

\address{Victoria University, PO Box 14428, Melbourne, Victoria, Australia}

\email{\url{Alasdair.McAndrew@vu.edu.au}}



\thanks{A version of this paper has been submitted to the International Journal of
  Computer Mathematics.\\}

\maketitle

\begin{abstract}
  Runge-Kutta formulas are some of the workhorses of numerical solving of
  differential equations.  However, they are extremely difficult to generate; the
  algebra involved can be very complicated indeed.  It is now standard, following
  with work of Butcher~\cite{butc87,butc08,lamb91} in the 1960's and 70's, to use the
  theory of trees to simplify the algebra.  More recently, however, several
  authors~\cite{gand99,grun95,fame04,sofr94} have shown that it is quite feasible to
  use a computer algebra system to generate Runge-Kutta formulas.  This article shows
  that, instead of using commercial systems as has been done previously, open-source
  systems can be used with equal effect.  This has the advantage that anybody can
  experiment with the code.\medskip
\end{abstract}

\section{Introduction}
\label{sec:intro}

We are concerned with finding the solution to the initial value problem
\[
\frac{dy}{dx}=f(x,y),\quad y(x_0)=y_0
\]
where the function $f$, and the initial values $(x_0,y_0)$ are given.  A
\emph{numerical solution} consists of a sequence of ordered pairs $(x_k,y_k)$, where
$y_k$ is an approximation to the exact value $y(x_k)$.  One way is to use the Taylor
expansion of $y(x)$, using
\begin{align*}
  y'&=f(x,y)\\
  y''&=f_x+f_y\frac{dy}{dx}=f_x+f_yf\\
  y'''&=(f_x+f_yf)_x+(f_x+f_yf)_y\frac{dy}{dx}\\
      &=f_{xx}+f_{yx}f+f_yf_x+(f_{xy}+f_{yy}f+f_yf_y)f\\
      &=f_{xx}+2f_{xy}f+f_xf_y+f_{yy}(f)^2+(f_y)^2f
\end{align*}
and then for a suitably small value of $h$, given $(x_k,y_k)$ and with
$x_{k+1}=x_k+h$ compute an approximation to $y_{k+1}\approx y(x_k+h)$ to the exact
value $y(x_{k+1})$.

However, this requires the derivatives of $f$, which in many cases may have to be
computed using an approximation, thus introducing a new source of errors.

The insight of Runge\footnote{Carl David Tolm\'e Runge, 1856--1927} and of
Kutta\footnote{Martin Wilhelm Kutta, 1867--1944}, was to realize that as the first
derivative of $y$ was equal to $f$, so other derivatives could be computed by
judicious nesting.  

For example, suppose we truncate the expansion of the Taylor series expansion of $f$
after the first derivative:
\begin{equation}
  \label{eq:tayl_f}
  f(x+h,y+k)\approx f+hf_x+kf_y.
\end{equation}
Also, truncate the Taylor series for $y$ after the second derivative:
\begin{equation}
  \label{eq:tayl_y}
  y(x+h)\approx y+hf+\frac{h^2}{2}(f_x+f_yf).
\end{equation}
Note that the expression in parentheses on the right of~\eqref{eq:tayl_y} is very close
to that on the right of~\eqref{eq:tayl_f}, excepting a term of $f$.  But this can be
inserted simply by writing~\eqref{eq:tayl_y} as
\begin{equation}
  \label{eq:tayl_y2}
  y(x+h)\approx y+\frac{h}{2}f+\frac{h}{2}(f+hf_x+hf_yf).
\end{equation}
Comparing the final term with~\eqref{eq:tayl_f} we can write
\begin{equation}
  \label{eq:tayl_f2}
  y(x+h)\approx y+\frac{h}{2}f+\frac{h}{2}f(x+h,y+hf).
\end{equation}
This can be written as a sequence of steps, starting with $y_n\approx y(x_n)$ and
with $x_{n+1}=x_n+h$:
\begin{align*}
  k_1&=f(x,y)\\
  k_2&=f(x+h,y+hk_1)\\
  y_{n+1}&=y_n+\frac{h}{2}(k_1+k_2).
\end{align*}
This is an example of a second-order Runge-Kutta formula, and is equal to a
second-order Taylor approximation, but without computing any of the derivatives of
$f$.  In general, an $n$-th order Runge-Kutta formula has the form:
\begin{align*}
  k_1&=f(x_n,y_n)\\
  k_2&=f(x+c_2h,y_n+a_{21}hk_1)\\
  k_3&=f(x+c_3h,y_n+h(a_{31}k_1+a_{32}k_2))\\
  &\vdots\\
  k_m&=f(x+c_mh,y_n+h(a_{m1}k_1+a_{m2}k_2+\cdots+a_{m,m-1}k_{m-1}))\\
  \intertext{and then}
  y_{n+1}&=y_n+h(b_1k_1+b_2k_2+\cdots+b_mk_m).
\end{align*}

It is customary to write all the coefficients in a \emph{Butcher array}:
\[
\begin{array}[h]{c|ccccc}
  0&&&&&\\
  c_2&a_{21}&&&\\
  c_3&a_{31}&a_{32}&&\\
  \vdots&&&&\\
  c_m&a_{m1}&a_{m2}&\cdots&a_{m,m-1}\\
  \hline
  &b_1&b_2&\cdots&b_{m-1}&b_m
\end{array}
\]
These particular Runge-Kutta methods are called \emph{explicit} methods, where
at stage the value $k_i$ is explicitly defined in terms of previously computed
values.  The above second order method could be written as
\[
\setlength{\arraycolsep}{5pt}
\begin{array}[h]{c|cc}
  0&&\\
  1&1&\\
  \hline
  &\arraystrut\frac{1}{2}&\frac{1}{2}
\end{array}
\]
A very popular fourth order method (sometimes called ``\emph{the} Runge-Kutta
method'') is given by the array
\[
\setlength{\arraycolsep}{5pt}
\begin{array}[h]{c|cccc}
  0&&&&\\
  \arraystrut\frac{1}{2}&\frac{1}{2}&&&\\
  \arraystrut\frac{1}{2}&0&\frac{1}{2}&&\\
  \arraystrut 1&0&0&1&\\
  \hline
  \arraystrut&\frac{1}{6}&\frac{1}{3}&\frac{1}{3}&\frac{1}{6}
\end{array}
\]

\section{Use of Computer Algebra Systems: third order methods}
\label{sec:cas}

Before we launch into the use of a CAS, consider a third-order system:
\begin{align*}
  k_1&=f(x,y)\\
  k_2&=f(x+c_2h,y+a_{21}hk_1)\\
  k_3&=f(x+c_3h,y+a_{31}hk_1+a_{32}hk_2)\\
  \intertext{and}
  y_{n+1}&=y_n+h(b_1k_1+b_2k_2+b_3k_3)
\end{align*}
and which is to be equivalent to the third-order Taylor polynomial
\[
y(x_n+h)\approx
y+hf+\frac{h^2}{2}(f_x+f_yf)+\frac{h^3}{6}(f_{xx}+2f_{xy}f+f_xf_y+f_{yy}(f)^2+(f_y)^2f).
\]
In order to find appropriate coefficients, the expressions for each of the $k_i$
values need to be expanded up to and including the second derivatives; thus:
\begin{align*}
  k_1&=f(x,y)\\
  k_2&=f(x,y)+h(c_2f_x+a_{21}k_1f_y)+\frac{h^2}{2}(c_2^2f_{xx}+2c_2a_{21}k_1f_{xy}
  +(a_{32}k_1)^2f_{yy})\\
  k_3&=f(x,y)+h(c_3f_x+(a_{31}k_1+a_{32}k_2)f_y)\\
  &\qquad +\frac{h^2}{2}(c_2^2f_{xx}+2c_2(a_{31}k_1+a_{32}k_2)f_{xy}
  +(a_{31}k_1+a_{32}k_2)^2f_{yy})
\end{align*}
Note that $k_3$ above is expressed in terms of $k_2$; this means that the expression
for $k_2$ must be substituted into $k_3$ wherever it occurs, so that the final
expressions for each of the $k_i$ are written using only $f$ and its derivatives, and
$h$.

To equate the Taylor polynomial with the Runge-Kutta values for $y_{n+1}$, we must
have
\begin{multline*}
  y+hf+\frac{h^2}{2}(f_x+f_yf)+\frac{h^3}{6}(f_{xx}+2f_{xy}f+f_xf_y+f_{yy}(f)^2+(f_y)^2f)\\
  = y+h(b_1k_1+b_2k_2+b_3k_3)
\end{multline*}
or that
\begin{equation}
  b_1k_1+b_2k_2+b_3k_3
  =f+\frac{h}{2}(f_x+f_yf)+\frac{h^2}{6}(f_{xx}+2f_{xy}f+f_xf_y+f_{yy}(f)^2+(f_y)^2f).\label{eq:rk1}
\end{equation}
We thus need to find values of all the unknown coefficients (the $a$, $b$ and $c$
values), for which
\[
b_1k_1+b_2k_2+b_3k_3-f-\frac{h}{2}(f_x+f_yf)-\frac{h^2}{6}(f_{xx}+2f_{xy}f+f_xf_y+f_{yy}(f)^2+(f_y)^2f)=0.
\]
Collecting all the terms together, and working through all the algebra to expand
$k_2$ and $k_3$ fully, we end up with

\begin{equation}\label{eq:rk_expand}
\begin{aligned}
  \MoveEqLeft[4] \left(\frac{a_{32}^2b_3}{2}+a_{31}a_{32}b_3+\frac{a_{31}^2b_3}{2}+
    \frac{a_{21}^2b_2}{2}-\frac{1}{6}\right)f_{yy}h^2f^2
+\left(a_{21}a_{32}b_3-\frac{1}{6}\right)f_{y}^2h^2f\\
&+\left(a_{32}b_{3}c_{3}+a_{31}b_{3}c_{3}+a_{21}b_{2}c_{2}-\frac{1}{3}\right)f_{xy}h^2f\\
&+\left(a_{32}b_{3}+a_{31}b_{3}+a_{21}b_{2}-\frac{1}{2}\right)f_{y}hf
+\left(b_3+b_2+b_1-1\right)f\\
&+\left(a_{32}b_3c_2-\frac{1}{6}\right)f_xf_yh^2
+\left(\frac{b_3c_3^2}{2}+\frac{b_2c_2^2}{2}-\frac{1}{6}\right)f_{xx}h^2\\
&+ \left(b_3c_3+b_2c_2-\frac{1}{2}\right)f_xh\\
&=0.
\end{aligned}
\end{equation}

Once this has been done, the values we want are the solutions to the non-linear
equations:
\begin{align*}
a_{32}^2b_3+2a_{31}a_{32}b_3+a_{31}^2b_3+a_{21}^2b_2&=1/3\\
a_{21}a_{32}b_3&=1/6\\
a_{32}b_{3}c_{3}+a_{31}b_{3}c_{3}+a_{21}b_{2}c_{2}&=1/3\\
a_{32}b_{3}+a_{31}b_{3}+a_{21}b_{2}&=1/2\\
b_3+b_2+b_1&=1\\
a_{32}b_3c_2&=1/6\\
b_3c_3^2+b_2c_2^2&=1/3\\
b_3c_3+b_2c_2&=1/2
\end{align*}
It can be seen---even without attempting to solve these equations---that the algebra
involved is extremely involved, messy, and without any apparent order.  One of the
remarkable advances made by Butcher was to relate all these equations to the theory
of rooted trees, and hence bring some order to the apparent chaos.

Our approach, though, will be to simply create the equations from scratch, and solve
them, using the open-source computer algebra system Sage~\cite{sage} to perform the
computations.  However, much of the initial calculus computations will devolve to
Maxima~\cite{maxima}, which is the current descendant of the venerable system
Macsyma; and which has very powerful calculus and algebra functionality.  As Sage
includes Maxima within it, we can use Maxima initially to create the derivatives and
the functions, and use the algebraic power of Sage to solve the equations.  We will
present our work with monospaced input and typeset output, similar to the appearance
of using Sage in a browser-based ``notebook''~\cite{eroc10} and with Maxima cells..

To start, we need to create a formal function$f$ and its derivatives.  In Sage,
objects maintain their types, so in a variable assignment such as ``z =
\texttt{f.diff(x)}'', assuming \texttt{f} to be a function previously defined, the
result \texttt{z} will be either a Sage object or a Maxima object depending on the
type of \texttt{f}.  We thus start by introducing three Maxima variables:
\begin{Verbatim}
x = maxima('x')
y = maxima('y')
f = maxima('f')
\end{Verbatim}
Now each of these variables will automatically have access to the Maxima sub-system,
and so we can create the derivatives.  In order to prevent unnecessary high derivatives
of $y(x)$, we shall replace $y'$ with $f$ as soon as it appears.
\begin{Verbatim}
y.depends(x)
f.depends([x,y])
f1 = f.diff(x).subst("diff(y,x)=f")
f2 = f1.diff(x).subst("diff(y,x)=f")
f3 = f2.diff(x).subst("diff(y,x)=f")
\end{Verbatim}
Sage doesn't display the results of a variable assignment, but we can check out the
first two:
\begin{Verbatim}
f1,f2
\end{Verbatim}

\begin{gather*}
  f\,\left({{{\it \partial}}\over{{\it \partial}\,y}}\,f\right)+{{
      {\it \partial}}\over{{\it \partial}\,x}}\,f,\\
  f\,\left(f\,\left({{{\it \partial}^2}\over{{\it \partial}\,y^2}}\,f
    \right)+{{{\it \partial}^2}\over{{\it \partial}\,x\,{\it \partial}\,
        y}}\,f\right)+{{{\it \partial}}\over{{\it \partial}\,y}}\,f\,\left(f
    \,\left({{{\it \partial}}\over{{\it \partial}\,y}}\,f\right)+{{
        {\it \partial}}\over{{\it \partial}\,x}}\,f\right)+{{{\it \partial}^
      2}\over{{\it \partial}\,x^2}}\,f+f\,\left({{{\it \partial}^2}\over{
        {\it \partial}\,x\,{\it \partial}\,y}}\,f\right)
\end{gather*}

\N{}Next we make some substitutions for easier work later on, first introducing some
variables into the namespace.  Being based on Python, any variable must be named
before it can be used.  These variables will accrue their appropriate types later.

\begin{Verbatim}
var('h,F,Fx,Fy,Fxx,Fxy,Fyy,Fxxx,Fxxy,Fxyy,Fyyy,a21,a31,a32,\
  a41,a42,a43,b1,b2,b3,b4,c2,c3,c4,')
dsubs = " 'diff(f,x,3)=Fxxx, 'diff(f,x,2,y,1)=Fxxy,\
  'diff(f,x,1,y,2)=Fxyy,'diff(f,y,3)=Fyyy, 'diff(f,x,2)=Fxx,\
  'diff(f,y,2)=Fyy, 'diff(f,x,1,y,1)=Fxy,'diff(f,x,1)=Fx,\
  'diff(f,y,1)=Fy, f=F"
F1 = f1.subst(dsubs)
F2 = f2.subst(dsubs)
F3 = f3.subst(dsubs)
\end{Verbatim}
As before, their values can be checked:
\begin{Verbatim}
F1, F2, F3
\end{Verbatim}
\begin{gather*}
  {\it Fy}\,F+{\it Fx}\\
  F\,\left({\it Fyy}\,F+{\it Fxy}\right)+{\it Fy}\,\left({\it Fy}\,F+  {\it
      Fx}\right)+{\it Fxy}\,F+{\it Fxx}\\
  F\,\left(F\,\left({\it Fyyy}\,F+{\it Fxyy}\right)+{\it Fyy}\,\left(  {\it
        Fy}\,F+{\it Fx}\right)+{\it Fxyy}\,F+{\it Fxxy}\right)\\
  {}\qquad +{\it Fy}  \,\left(F\,\left({\it Fyy}\,F+{\it Fxy}\right)+{\it Fy}\,\left(  {\it
        Fy}\,F+{\it Fx}\right)+{\it Fxy}\,F+{\it Fxx}\right)\\
  {}\qquad +2\,\left(  {\it Fy}\,F+{\it Fx}\right)\,\left({\it Fyy}\,F+{\it Fxy}\right)+  {\it
    Fxy}\,\left({\it Fy}\,F+{\it Fx}\right)+F\,\left({\it Fxyy}\,F+  {\it
      Fxxy}\right)\\
  {}\qquad +{\it Fxxy}\,F+{\it Fxxx}
  \end{gather*}
Now we introduce the Taylor polynomial up to the third derivative (this is for a
third order method), and this corresponds to the right hand side of
equation~\eqref{eq:rk1}:
\begin{Verbatim}
T = F + h/2*F1 + h^2/6*F2; T
\end{Verbatim}
\[
{{h^2\,\left(F\,\left({\it Fyy}\,F+{\it Fxy}\right)+{\it Fy}\,  \left({\it Fy}\,F+{\it Fx}\right)+{\it Fxy}\,F+{\it Fxx}\right)  }\over{6}}+{{h\,\left({\it Fy}\,F+{\it Fx}\right)}\over{2}}+F
\]
In order to compute the $k_i$ values, we need a Taylor series expansion up to the
second derivative:
\[
f(x+a,y+b)=f(x,y)+af_x+bf_y+\frac{1}{2}\left(a^2f_{xx}+2abf_{xy}+b^2f_{yy}\right)
\]
where the subscripts represent the usual partial derivatives.  Given the above
substitutions, we will call this expansion $\mathrm{Tay}(a,b)$.  Since we are at the
moment dealing with Maxima variables, we will define Tay as a Maxima function:
\begin{Verbatim}
Tay = maxima.function('a,b','F+Fx*a+Fy*b+(Fxx*a^2+2*Fxy*a*b+Fyy*b^2)/2')
\end{Verbatim}
Since the $k_i$ values are nested, we don't want the powers of $h$ increasing: we are
only interested in coefficients for which the powers of $h$ are 2 or less.  Maxima
has a handy trick here:
\begin{Verbatim}
maxima("tellrat(h^3)")
maxima("algebraic:true")
\end{Verbatim}
This means that for every rational expansion, all powers of $h$ which are three or
more will be set equal to zero.  Now we can create the $k_i$ values:
\begin{Verbatim}
k1 = Tay(0, 0)
k2 = Tay(c2*h, a21*h*k1)
k3 = Tay(c3*h, h*a31*k1 + h*a32*k2)
\end{Verbatim}
(The expressions, certainly for $k_3$, are too long to display).  Now we can create
the left hand side of equation~\eqref{eq:rk1}:
\begin{Verbatim}
RK = b1*k1 + b2*k2 + b3*k3
\end{Verbatim}
The next step is to compute create the expression on the left hand side of
equation~\eqref{eq:rk_expand}; we can do this by collecting all the terms involving $F$
and its derivatives.  The Maxima command ``\texttt{collectterms}'' provides just this
functionality. 
\begin{Verbatim}
d = (T-RK).ratexpand().collectterms(Fyy,Fxy,Fxx,Fy,Fx,F,h)
\end{Verbatim}
This expression is too long to print, but we can extract the coefficients from it,
which are the equations we want:
\begin{Verbatim}
eqs = [xx.inpart(1) for xx in d.args()]
eqs
\end{Verbatim}

\begin{gather*}
  \left[-{{{\it a_{32}}^2\,{\it b_3}}\over{2}}-{\it a_{31}}\,{\it a_{32}}\,  {\it
      b_3}-{{{\it a_{31}}^2\,{\it b_3}}\over{2}}-{{{\it a_{21}}^2\,  {\it
          b_2}}\over{2}}+{{1}\over{6}},\quad {{1}\over{6}}-{\it a_{21}}\,{\it
 a_{32}}\,{\it b_3},\right.\\
{}\quad -{\it a_{32}}\,{\it b_3}\,{\it c_3}-{\it a_{31}}\,{\it b_3}\,  {\it c_3}-{\it
  a_{21}}\,{\it b_2}\,{\it c_2}+{{1}\over{3}},\quad -{\it a_{32}}\,{\it b_3}-{\it
  a_{31}}\,{\it b_3}-{\it a_{21}}\,  {\it b_2}+{{1}\over{2}},\\
{}\quad \left.-{\it b_3}-{\it b_2}-{\it b_1}+1, {{1}\over{6}}-{\it a_{32}}\,{\it b_3}\,{\it c_2}, -{{{\it b_3}\,{\it c_3}^2}\over{2}}-{{{\it b_2}\,{\it c_2}^2}\over{  2}}+{{1}\over{6}}, -{\it b_3}\,{\it c_3}-{\it b_2}\,{\it c_2}+{{1}\over{2}}\right]
\end{gather*}
These equations are all Maxima objects, but to access the algebraic power of Sage,
they need to be lifted out of the Maxima sub-system.  We will use the Sage
``\texttt{repr}'' command, which produces a string representation of an object, and
we will evaluate those strings into expressions.  

\begin{Verbatim}
eqs2 = [sage_eval(repr(xx),locals=locals()) for xx in eqs]
\end{Verbatim}

\N{}These equations can be solved in terms of $c_2$ and $c_3$:
\begin{Verbatim}
sols = solve(eqs2,[a21,a31,a32,b1,b2,b3])
sols
\end{Verbatim}

\begin{gather*}
\left[\left[a_{21} = c_{2},\quad a_{31} = \frac{3 \, {\left(c_{2}^{2} - c_{2}\right)}
      c_{3} + c_{3}^{2}}{3 \, c_{2}^{2} - 2 \, c_{2}},\quad a_{32} = \frac{c_{2}
      c_{3} - c_{3}^{2}}{3 \, c_{2}^{2} - 2 \, c_{2}},\right.\right.\\ 
    {}\quad\left.\left. b_{1} = \frac{3 \, {\left(2 \, c_{2} - 1\right)} c_{3} - 3 \, c_{2} + 2}{6 \,
      c_{2} c_{3}}, \quad b_{2} = -\frac{3 \, c_{3} - 2}{6 \, {\left(c_{2}^{2} -
          c_{2} c_{3}\right)}},\quad  b_{3} = \frac{3 \, c_{2} - 2}{6 \, {\left(c_{2} c_{3} - c_{3}^{2}\right)}}\right]\right]
\end{gather*}
These are standard expressions, and be found, for example, in Butcher~\cite{butc08}.
In general, Sage solutions are given as a list of lists.  In our case there is only
one solution, which can be isolated with
\begin{Verbatim}
sols = sols[0]
\end{Verbatim}
given that in Sage lists are indexed starting at zero.  Note that the first item
tells us that $a_{21}=c_2$.  Adding the next two items shows a similar relation for
$c_3$:
\begin{Verbatim}
(sols[1] + sols[2]).simplify_rational()
\end{Verbatim}
\[
a_{31} + a_{32} = c_{3}
\]
It can in fact be shown that for any Runge-Kutta method, each $c_k$ value is equal
to the sum of the corresponding $a_{ki}$ values:
\[
c_k=a_{k1}+a_{k2}+\cdots a_{k,k-1}.
\]
This is known as the \emph{row-sum condition}, and may be assumed for any computation
with Runge-Kutta coefficients.

We can now find particular solutions by substituting values for $c_2$ and $c_3$ (such
that the denominators are all non-zero, which means that $c_2$ and $c_3$ must be
different), for example:
\begin{Verbatim}
[xx.subs(c2=-1,c3=1) for xx in sols]
\end{Verbatim}
\[
\left[a_{21} = \left(-1\right), a_{31} = \left(\frac{7}{5}\right), a_{32} = \left(-\frac{2}{5}\right), b_{1} = \left(\frac{2}{3}\right), b_{2} = \left(-\frac{1}{12}\right), b_{3} = \left(\frac{5}{12}\right)\right]
\]
Two other substitutions are:
\begin{Verbatim}
[xx.subs(c2=1/2,c3=1) for xx in sols]
\end{Verbatim}
\[
\left[a_{21} = \left(\frac{1}{2}\right), a_{31} = \left(-1\right), a_{32} = 2, b_{1} = \left(\frac{1}{6}\right), b_{2} = \left(\frac{2}{3}\right), b_{3} = \left(\frac{1}{6}\right)\right]
\]
\begin{Verbatim}
[xx.subs(c2=1/2,c3=1) for xx in sols]
\end{Verbatim}
\[
\left[a_{21} = \left(\frac{1}{3}\right), a_{31} = 0, a_{32} =
  \left(\frac{2}{3}\right), b_{1} = \left(\frac{1}{4}\right), b_{2} = 0, b_{3} =
  \left(\frac{3}{4}\right)\right]
\]
These last two have Butcher arrays
\[
\setlength{\arraycolsep}{5pt}
\begin{array}[h]{c|ccc}
  0&&&\\
  \arraystrut\frac{1}{2}&\frac{1}{2}&&\\
  \arraystrut 1&-1&2&\\
  \hline
  \arraystrut&\frac{1}{6}&\frac{2}{3}&\frac{1}{6}
\end{array},\qquad
\setlength{\arraycolsep}{5pt}
\begin{array}[h]{c|ccc}
  0&&&\\
  \arraystrut\frac{1}{3}&\frac{1}{3}&&\\
  \arraystrut\frac{2}{3}&0&\frac{2}{3}&\\
  \hline
  \arraystrut&\frac{1}{3}&0&\frac{3}{4}
\end{array}
\]
and are known as \emph{Kutta's third-order method} and \emph{Heun's third-order
  method} respectively.

\section{Fourth-order methods}
\label{sec:fourth}

Having set up the ground work, fourth order methods can be found similarly; with
suitable changes to some of the entries to allow for the higher order.  With each of
\texttt{f1, f2, f3, F1, F2, F3} as before, the commands (given with no outputs) will
be:
\begin{Verbatim}
maxima("tellrat(h^4)")
T = F + h/2*F1 + h^2/6*F2 + h^3/24*F3
Tay = maxima.function('a,b','F+Fx*a+Fy*b+(Fxx*a^2+2*Fxy*a*b+Fyy*b^2)/2\
     +(Fxxx*a^3+3*Fxxy*a^2*b+3*Fxyy*a*b^2+Fyyy*b^3)/6')
k1 = Tay(0, 0)
k2 = Tay(c2*h, a21*h*k1)
k3 = Tay(c3*h, h*a31*k1 + h*a32*k2)
k4 = Tay(c4*h, h*a41*k1 + h*a42*k2 + h*a43*k3)

RK = b1*k1 + b2*k2 + b3*k3 + b4*k4
d = (T-RK).ratexpand()\
    .collectterms(Fyyy,Fxyy,Fxxy,Fxxx,Fyy,Fxy,Fxx,Fy,Fx,F,h)
eqs = [xx.inpart(1) for xx in d.args()]
\end{Verbatim}
At this stage, with no simplification, we will have a list of 19 equations.  Before
the business of simplification, we first introduce the row-sum condition:
\begin{Verbatim}
eqs2 = [xx.subst("a21=c2,a31=c3-a32,a41=c4-a42-a43").expand()\
    .collectterms(b1,b2,b3,b4) for xx in eqs]
\end{Verbatim}
and transform the set of equations out of Maxima and into Sage:
\begin{Verbatim}
eqss = [sage_eval(repr(xx),locals=locals()) for xx in eqs2]
\end{Verbatim}
Now we can create a polynomial ring in which all the computations will be done, and
in the polynomial ring compute the reduced basis of the ideal generated by the
equations:
\begin{Verbatim}
R = PolynomialRing(QQ,'a21,a31,a32,a41,a42,a43,b1,b2,b3,b4,c2,c3,c4',order='lex')
Id = R.ideal(eqss)
ib = Id.interreduced_basis()
\end{Verbatim}
This is now a list of only eight equations, instead of the 19 from earlier:
\begin{gather*}
a_{32} a_{43} b_{4} c_{2} - \frac{1}{24}\\
a_{32} b_{3} c_{2} + a_{42} b_{4} c_{2} + a_{43} b_{4} c_{3} -
\frac{1}{6}\\
a_{42} b_{4} c_{2} c_{3} -  a_{42} b_{4} c_{2} c_{4} + a_{43} b_{4}
c_{3}^{2} -  a_{43} b_{4} c_{3} c_{4} - \frac{1}{6} c_{3} + \frac{1}{8}\\
a_{43} b_{4} c_{2} c_{3} -  a_{43} b_{4} c_{3}^{2} - \frac{1}{6} c_{2} +
\frac{1}{12}\\
b_{1} + b_{2} + b_{3} + b_{4} - 1\\
b_{2} c_{2} + b_{3} c_{3} + b_{4} c_{4} - \frac{1}{2}\\
b_{3} c_{2} c_{3} -  b_{3} c_{3}^{2} + b_{4} c_{2} c_{4} -  b_{4}
c_{4}^{2} - \frac{1}{2} c_{2} + \frac{1}{3}\\
b_{4} c_{2} c_{3} c_{4} -  b_{4} c_{2} c_{4}^{2} -  b_{4} c_{3}
c_{4}^{2} + b_{4} c_{4}^{3} - \frac{1}{2} c_{2} c_{3} + \frac{1}{3}
c_{2} + \frac{1}{3} c_{3} - \frac{1}{4}\\  
\end{gather*}
And this new system of equations can be easily solved in terms of $c_i$.  We first
note that $c_4=1$; although this can be shown analytically, we can easily demonstrate
it using our reduced Gr\"obner basis:
\begin{Verbatim}
(ib*R).reduce(c4-1)
\end{Verbatim}
\[
0
\]
The $b_i$ values can now be obtained by noting that the last four equations are
linear in $b_i$:
\begin{Verbatim}
bsols = solve([ib[i].subs(c4=1) for i in [4,5,6,7]],[b1,b2,b3,b4],\
    solution_dict=True)[0]
bsols = {xx: factor(yy) for xx, yy in bsols.items()}
\end{Verbatim}
\begin{multline*}
  \left\{b_{2} : \frac{2 \, c_{3} - 1}{12 \, {\left(c_{2} - c_{3}\right)}
      {\left(c_{2} - 1\right)} c_{2}}, b_{1} : \frac{6 \, c_{2} c_{3} - 2 \, c_{2} -
      2 \, c_{3} + 1}{12 \, c_{2} c_{3}},\right.\\
    \left. b_{4} : \frac{6 \, c_{2} c_{3} - 4 \, c_{2}
      - 4 \, c_{3} + 3}{12 \, {\left(c_{2} - 1\right)} {\left(c_{3} - 1\right)}},
    b_{3} : -\frac{2 \, c_{2} - 1}{12 \, {\left(c_{2} - c_{3}\right)} {\left(c_{3} -
          1\right)} c_{3}}\right\}
\end{multline*}
These are standard results~\cite{butc08}.  From these we can use the first equations
to compute the $a_{ij}$ values:
\begin{Verbatim}
asols = solve([SR(ib[i]).subs(c4=1).subs(bsols) for i in [1,2,3]],\
    [a32,a42,a43],solution_dict=True)[0]
asols = {xx: factor(yy) for xx, yy in asols.items()}
\end{Verbatim}

\begin{multline*}
  \left\{a_{43} : \frac{{\left(2 \, c_{2} - 1\right)} {\left(c_{2} - 1\right)}
      {\left(c_{3} - 1\right)}}{{\left(6 \, c_{2} c_{3} - 4 \, c_{2} - 4 \, c_{3} +
          3\right)} {\left(c_{2} - c_{3}\right)} c_{3}},\right.\\
  \left. a_{42} : -\frac{{\left(4 \,
          c_{3}^{2} - c_{2} - 5 \, c_{3} + 2\right)} {\left(c_{2} - 1\right)}}{2 \,
      {\left(6 \, c_{2} c_{3} - 4 \, c_{2} - 4 \, c_{3} + 3\right)} {\left(c_{2} -
          c_{3}\right)} c_{2}}, a_{32} : \frac{{\left(c_{2} - c_{3}\right)} c_{3}}{2
      \, {\left(2 \, c_{2} - 1\right)} c_{2}}\right\}
\end{multline*}
These values satisfy the first equation:
\begin{Verbatim}
(asols[a32]*asols[a43]*bsols[b4]*c2).rational_simplify()
\end{Verbatim}
\[
\frac{1}{24}
\]
Note that because of the factor $c_2-c_3$ in the denominators of some of these
expressions, we can't substitute equal values for $c_2$ and $c_3$.  In order to
develop a fourth-order method for which $c_2=c_3$ we need to go back a few steps:
\begin{Verbatim}
var('u')
eqss = [sage_eval(repr(xx),locals=locals()).subs(c2=u,c3=u,c4=1) for xx in eqs2]
R.<a21,a31,a32,a41,a42,a43,b1,b2,b3,b4,u> = PolynomialRing(QQ)
Id = R.ideal(eqss)
ib = Id.interreduced_basis(); ib
\end{Verbatim}
\begin{gather*}
\left[a_{32} a_{43} - \frac{1}{2}, a_{32} b_{3} - \frac{1}{6}, a_{42} +
a_{43} - 1, b_{1} - \frac{1}{6}, b_{2} + b_{3} - \frac{2}{3}, b_{4} -
\frac{1}{6}, u - \frac{1}{2}\right]  
\end{gather*}
These can be solved to produce:
\[
u=\frac{1}{2},\quad b_1=\frac{1}{6},\quad b_2=r_1,\quad b_3=\frac{2}{3}-r_1,\quad b_4=\frac{1}{6},\quad
a_{32}=\frac{1}{2(2-3r_1)},\quad a_{42}=3r_1-1,\quad a_{43}=2-3r_1.
\]
The Butcher array corresponding to this is
\[
\setlength{\arraycolsep}{5pt}
\begin{array}[h]{c|cccc}
  0&&&&\\
  \arraystrut \frac{1}{2}&\frac{1}{2}&&&\\
  \arraystrut \frac{1}{2}&0&\frac{1}{2}&&\\
  \arraystrut 1&0&3r_1-1&2-3r_1&\\
  \hline
  \arraystrut&\frac{1}{6}&r_1&\frac{2}{3}-r_1&\frac{1}{6}
\end{array}
\]
Putting $r_1=1/3$ produces the classic Runge-Kutta fourth order method.

\section{Use of autonomy}
\label{sec:auton}

Much of the computations in the previous sections can be simplified by noting that we
don't in fact need to include the $c_i$ values in any equations list, as their values
can be determined from the $a_{ij}$ values.  Since $c_i$ only appear in the
computation of $k_i$, all the computations can be simplified by considering only
differential equations for the form
\[
y'=f(y)
\]
in which $f$ does not depend explicitly on $x$; such differential equations are said
to be \emph{autonomous}.  This leads to greatly simplified forms for the higher
derivatives of $f$:
\begin{Verbatim}
f.depends(y)
\end{Verbatim}
Create $f_i$ and $F_i$ as before, but note the values of $F_i$:
\begin{Verbatim}
F1, F2, F3
\end{Verbatim}
\[
\left({\it Fy}\,F,\quad
  {\it Fyy}\,F^2+{\it Fy}^2\,F,\quad
  {\it Fyyy}\,F^3+4\,{\it
Fy}\,{\it Fyy}\,F^2+{\it Fy}^3\,F\right)
\]
Since we need not consider any partial derivatives of $f$ which include $x$, the
following commands can be used:
\begin{Verbatim}
maxima("tellrat(h^4)")
T = F + h/2*F1 + h^2/6*F2 + h^3/24*F3
Tay = maxima.function('a,b','F+Fy*b+Fyy*b^2/2+Fyyy*b^3/6')
k1 = Tay(0, 0)
k2 = Tay(0, a21*h*k1)
k3 = Tay(0, h*a31*k1 + h*a32*k2)
k4 = Tay(0, h*a41*k1 + h*a42*k2 + h*a43*k3)
\end{Verbatim}
Note that we need to include an extra dummy variable in the ``Tay'' function; the
systems in their current forms prevent a Maxima function of one variable being used
in this way.  The next few commands are similar to those above.
\begin{Verbatim}
RK = b1*k1 + b2*k2 + b3*k3 + b4*k4
d = (T-RK).ratexpand().collectterms(Fyyy,Fyy,Fy,F,h)
eqs = [xx.inpart(1) for xx in d.args()]
eqss = [sage_eval(repr(xx),locals=locals()) for xx in eqs]
\end{Verbatim}
The equations are still quite long and complicated; as before they can be simplified
by introducing the row-sum conditions, and putting $c_4=1$:
\begin{Verbatim}
eqs2 = [xx.subs(a31=c3-a32,a41=1-a42-a43).expand() for xx in eqss]
\end{Verbatim}
to produce:
\begin{gather*}
-\frac{1}{6} \, b_{2} c_{2}^{3} - \frac{1}{6} \, b_{3} c_{3}^{3} -
\frac{1}{6} \, b_{4} c_{4}^{3} + \frac{1}{24}\\
-\frac{1}{2} \, a_{32} b_{3} c_{2}^{2} - \frac{1}{2} \, a_{42} b_{4}
c_{2}^{2} - a_{32} b_{3} c_{2} c_{3} - \frac{1}{2} \, a_{43} b_{4}
c_{3}^{2} - a_{42} b_{4} c_{2} c_{4} - a_{43} b_{4} c_{3} c_{4} +
\frac{1}{6}\\
-\frac{1}{2} \, b_{2} c_{2}^{2} - \frac{1}{2} \, b_{3} c_{3}^{2} -
\frac{1}{2} \, b_{4} c_{4}^{2} + \frac{1}{6}\\
-a_{32} a_{43} b_{4} c_{2} + \frac{1}{24}\\
-a_{32} b_{3} c_{2} - a_{42} b_{4} c_{2} - a_{43} b_{4} c_{3} +
\frac{1}{6}\\
-b_{2} c_{2} - b_{3} c_{3} - b_{4} c_{4} + \frac{1}{2}\\
-b_{1} - b_{2} - b_{3} - b_{4} + 1
\end{gather*}
This equations are not quite the same as those from the previous section, but they
can be solved similarly:
\begin{Verbatim}
bsols = solve([eqs2[i] for i in [0,2,5,6]],[b1,b2,b3,b4],\
    solution_dict=True)[0]
asols = solve([eqs2[i].subs(bsols) for i in [1,3,4]],[a32,a42,a43],\
    solution_dict=True)[1]
\end{Verbatim}
to produce the same results as before.  

For this autonomous approach, there has been no need to simplify a large set of
nineteen equations to a smaller set by involving the machinery of Gr\"obner bases;
the equation set was optimally small at the start.

To solve these equations with $c_2=c_3$, we need to be a bit careful; the attempt
\begin{Verbatim}
var('u')
eqs2 = [xx.subs(a21=u,a31=u-a32,a41=1-a42-a43).expand() for xx in eqss]
bsols = solve([eqs2[i] for i in [0,2,5,6]],[b1,b2,b3,b4],solution_dict=True)
\end{Verbatim}
will not work: as $b_2$ and $b_3$ have the same coefficients in all the equations,
the determinant of the matrix of coefficients is zero.  So we leave $b_2$ out:
\begin{Verbatim}
bsols = solve([eqs2[i] for i in [0,2,6]],[b1,b3,b4],solution_dict=True)[0]
asols = solve([eqs2[i].subs(bsols) for i in [1,3,4]],[a32,a42,a43],\
    solution_dict=True)[0]
\end{Verbatim}
The results will be expressed in terms of the parameters $b_2$ and $u$.  Substituting
$b_2=1/3$ and $u=1/2$ will produce the standard fourth-order method.

\section{Embedded formulas}
\label{sec:embed}

Many applications now use \emph{embedded Runge-Kutta methods}, in which two methods
share the same coefficients.  Generally the order of the methods differs by one, so
we might have a fifth order method, from which the coefficients can be used to build
a fourth order method.  Then the differences between the results of these methods can
be used to adjust the step size $h$ for the next iteration.

Methods of order $4(3)$ and the theory behind them are well
known~\cite{hair87,butc08}, but we can show that it is very easy to construct such a
method.  Starting with Kutta's 3/8 method, we need to find coefficients
$\hat{b}_1,\hat{b}_2,\hat{b}_3,\hat{b}_4,\hat{b}_5$ so that with the extra stage
\begin{align*}
  k_5&=f(x_n+c_5h, h(b_1k_1+b_2k_2+b_3k_3+b_4k_4))\\
  \intertext{the value of $y_{n+1}$ obtained with}
  y_{n+1}&=y_n+h(\hat{b}_1k_1+\hat{b}_2k_2+\hat{b}_3k_3+\hat{b}_4k_4+\hat{b}_5k_5)
\end{align*}
will be accurate to order three.

This is easily done, assuming the autonomous approach.  We enter the fourth order
values, and for simplicity we use $s_i$ in place of $\hat{b}_i$.
\begin{Verbatim}
a21,a31,a32,a41,a42,a43 = 1/3,-1/3,1,1,-1,1
b1,b2,b3,b4 = 1/8,3/8,3/8,1/8
var('s1,s2,s3,s4,s5')
\end{Verbatim}
We set up the third order conditions as previously, but this time with five stages:
\begin{Verbatim}
T = F + h/2*F1 + h^2/6*F2
maxima("tellrat(h^3)")
maxima("algebraic:true")
Tay = maxima.function('a,b','F+Fy*b+Fyy*b^2/2')

k1 = Tay(0, 0)
k2 = Tay(0, h*a21*k1)
k3 = Tay(0, h*a31*k1 + h*a32*k2)
k4 = Tay(0, h*a41*k1 + h*a42*k2 + h*a43*k3)
k5 = Tay(0, h*b1*k1 + h*b2*k2 + h*b3*k3 + h*b4*k4)
RK = s1*k1 + s2*k2 + s3*k3 + s4*k4 + s5*k5
\end{Verbatim}
Now we extract the coefficients of \texttt{T-RK} as equations to be solved.
\begin{Verbatim}
d = (T-RK).ratexpand().collectterms(Fyy,Fy,F,h)
eqs = [xx.subst('h=1,F=1,Fy=1,Fyy=1') for xx in d.args()]
eqs2 = [sage_eval(repr(xx),locals=locals()) for xx in eqs]
\end{Verbatim}
These equations are easily solved:
\begin{Verbatim}
solve(eqs2,[s1,s2,s3,s4,s5])
\end{Verbatim}
\[
\left[\left[s_{1} = -\frac{1}{4} \, r_{1} + \frac{1}{8}, s_{2} =
\frac{3}{4} \, r_{1} + \frac{3}{8}, s_{3} = -\frac{3}{4} \, r_{1} +
\frac{3}{8}, s_{4} = -\frac{3}{4} \, r_{1} + \frac{1}{8}, s_{5} =
r_{1}\right]\right]
\]
and the extra parameter can be set to any value we like, for example $r_1=1$.
\begin{Verbatim}
[xx.subs(r1=1) for xx in sols[0]]
\end{Verbatim}
\[
\left[s_{1} = \left(-\frac{1}{8}\right), s_{2} =
\left(\frac{9}{8}\right), s_{3} = \left(-\frac{3}{8}\right), s_{4} =
\left(-\frac{5}{8}\right), s_{5} = 1\right]
\]
This leads to a Butcher array
\[
\setlength{\arraycolsep}{5pt}
\begin{array}[h]{c|rrrrr}
  0&&&&&\\
  \arraystrut \frac{1}{3}&\frac{1}{3}&&&&\\
  \arraystrut \frac{2}{3}&-\frac{1}{3}&1&&\\
  \arraystrut 1&1&-1&1&\\
  \hline
  \arraystrut&\frac{1}{8}&\frac{3}{8}&\frac{3}{8}&\frac{1}{8}&\\
  \arraystrut&-\frac{1}{8}&\frac{9}{8}&-\frac{3}{8}&-\frac{5}{8}&1
\end{array}
\]
for an embedded $4(3)$ method.  If we choose the parameter $r_1$ so that $s_4=0$;
that is $r_1=1/6$, we obtain the values
\[
\left[s_{1} = \left(\frac{1}{12}\right), s_{2} =
\left(\frac{1}{2}\right), s_{3} = \left(\frac{1}{4}\right), s_{4} = 0,
s_{5} = \left(\frac{1}{6}\right)\right]
\]

\section{Conclusions}
\label{sec:conc}

The literature on Runge-Kutta methods and associated mathematics is vast,
Butcher~\cite{butc08} lists many hundreds of references.  However, much of this
material is directed at the specialist researcher.  The various articles which use
computer algebra systems in an attempt to sidestep the specialist material have
tended to use commercial systems, which puts the material out of bounds for people
who don't use (or can't afford) those systems.  As long ago as 1993, Joachim
Neub\"user, the creator of the GAP package for group theory (and which is part of
Sage)~\cite{neub95} deplored the fact that mathematical theorems are open to
everybody to use, but mathematics using a computer system was not, unless that system
was open-source.  Our article has attempted to allow the general non-specialist
reader to experiment and explore some of the basic properties of Runge-Kutta
methods, using only open-source systems.

\bibliographystyle{abbrv} \bibliography{refs}

\end{document}